%% file: main.tex
\documentclass[conf]{new-aiaa}
\usepackage[utf8]{inputenc}
\usepackage{booktabs} 
\usepackage{array}    
\usepackage{multirow} 
\usepackage{graphicx}
\usepackage{amsmath}
\usepackage[version=4]{mhchem}
\usepackage{siunitx}
\usepackage{longtable,tabularx}
\usepackage{placeins} 
\usepackage{float}
\title{Quantum-Assisted Space Logistics Mission Planning}

\author{Amiratabak Bahengam\footnote{Ph.D. Student, Department of Systems and Enterprises, 1 Castle Point Terrace, AIAA Student Member.}}
\affil{Stevens Institute of Technology, Hoboken, NJ, 07030}

\author{Mohammad-Ali Miri\footnote{Scientific Advisor, 5 Marine View Plaza.}, R. Joseph Rupert\footnote{Embedded Software Engineer, 5 Marine View Plaza.}, Wesley Dyk\footnote{Senior Quantum Solutions Architect, 5 Marine View Plaza.}} 
\affil{Quantum Computing Inc (QCi), Hoboken, NJ, 07030}
\author{Hao Chen\footnote{Assistant Professor, Department of Systems and Enterprises, 1 Castle Point Terrace, AIAA Member.}}
\affil{Stevens Institute of Technology, Hoboken, NJ, 07030}

\begin{document}

\maketitle

\begin{abstract}
Quantum computing provides a novel approach to addressing conventionally intractable issues in large-scale optimization. Space logistics missions require the efficient routing of payloads, spacecraft, and resources across complex networks, often resulting in an exponential growth of the solution space that classical methods cannot efficiently solve. This paper leverages entropy quantum computing to model and solve the space logistics problem as a time-dependent multicommodity network flow, enabling the exploration of large solution spaces. The findings highlight quantum computing's potential to address complex aerospace logistics, demonstrating its suitability for complex interplanetary mission planning.
\end{abstract}

\section{Nomenclature}

{\renewcommand\arraystretch{1.0}
\noindent\begin{longtable*}{@{}l @{\quad=\quad} l@{}}
$G$ & directed graph \\
$N$ & set of depots \\
$E$ & set of directed routes \\
$t$ & time step index \\
$k$ & index for commodities \\
$d_{ij}$ & distance of route $(i, j)$ \\
$c_{ij}$ & cost efficient per vehicle for arc \( (i,j) \) \\
$W$ & capacity of each vehicle \\
$L_k$ & load size of commodity $k$ \\
$d_{ki}^t$ & time-dependent supply (+) or demand (-) of commodity $k$ at depot $i$ at time $t$ \\
$x_{ij}^{k,t}$ & flow of commodity $k$ on route $(i, j)$ at time $t$ \\
$z_{ij}^t$ & number of vehicles on route $(i, j)$ at time $t$ \\
$\mathcal{J}$ & objective function
\end{longtable*}}

\section{Introduction}
\lettrine{Q}uantum computing is an emerging field that leverages the principles of quantum mechanics to revolutionize the way we solve computational problems. Unlike classical computing, which uses bits to encode information as either 0 or 1, quantum computing uses quantum bits, or qubits \cite{schumacher_quantum_1995}. Qubits can exist in a state of 0, 1, or any quantum superposition of these states, enabling quantum systems to process a vast amount of information simultaneously. This capability, combined with quantum entanglement, allows quantum computers to explore multiple solutions to a problem at once, offering a significant advantage for specific types of computations. In this study, we explore what advantages quantum computing may provide for addressing important issuies in the field of space logistics.

As space exploration attracts more and more attention from various stakeholders, there is a strong interest in establishing long-term human presence in space. This is different from previous mission planning demands, where space missions are typically designed independently. For example, in Apollo missions, we conduct a carry-long strategy to carry everything we need all-at-once for each lunar mission; in International Space Station supply missions, we conduct a resupply strategy to deliver payload upon demands each time. The development of space infrastructures, such as in situ resource utilization (ISRU), power system, and the lunar foundation surface habitat specified in Artemis missions \cite{united_states_national_aeronautics_and_space_administration_nasas_2021}, boost a demand to consider mission interdependence in space mission design. These infrastructures need to be deployed at the early stage of the mission, generating high mission costs; We anticipate the high initial deployment cost can be paid off in later stages of the mission by utilizing the infrastructures to better support the exploration (e.g., resource generation from ISRU). The diverse infrastructures for distinct mission objectives and long time horizon make the logistics planning even more challenging.

Current space logistics models have attempted to address these challenges. Chen et. al proposed a piecewise linear approximation method to achieve a linearized approximation solution for integration design of spacecraft, ISRU, and logistics scheduling \cite{chen_integrated_2018}. Takubo et. al developed a reinforcement learning-based framework to resolve the computation challenge by relying on intelligent agents to explore the problem solution space \cite{takubo_hierarchical_2022}. Multiple methods have been proposed to decompose the problem \cite{isaji_multidisciplinary_2022}, aggregate formulation constraints and variables to reduce its complexity \cite{chen_multifidelity_2021}, restructure the problem space through a hierarchical framework \cite{takubo_hierarchical_2022}, or consider aggregated vehicle concepts \cite{downs_spaceflight_2023}.

This paper will develop a quantum computing-based framework to resolve the computation challenges of long-term multi-mission space campaign planning. We discretize the problem space by proposing a Hamiltonian-based system state representation, which is a ubiquitous concept for physical systems. Then, we can perform adiabatic quantum optimization (AQO) to identify the optimal system Hamiltonian that optimizes the logistics planning. We will determine a few key mission planning decisions as Hamiltonian states. Other details logistics operation can be generated following the optimized planning decisions using network-based space logistics model, making it a quantum-assisted logistics mission planning. This paper will be one of the first attempts to bring quantum computing in solving large-scale space logistics planning problems. 

The first hybrid quantum-classical algorithm is the variational quantum eigensolver (VQE) \cite{tilly_variational_2022} to identify the ground state of chemical molecules. For classical computers, to represent the occupation state of $N$ molecular orbits, a $2^N$-dimension linear space is needed, whereas only $N$ qubits are required for a quantum computer to represent this design space. This is because of the principle of quantum superposition, where qubits can exist in states and at the same time.

In this paper, we present an approach to solve the time-dependent multicommodity network flow (MCNF) problem in the context of space logistics by leveraging quantum computing. Our method formulates the problem as a directed graph with temporal dynamics, transforming it into a Hamiltonian suitable for quantum optimization. By utilizing quantum mechanical principles, we efficiently explore the solution space, addressing the complexity of large-scale, resource-constrained logistics problems.

The remainder of this paper is structured as follows: Section \ref{sec: Methodology} details the methodology and implementation, including the formulation of the space logistics problem, its representation as a quantum-computable model, and the process of encoding and solving it using quantum computing techniques. In Section \ref{sec: Results and Analysis}, we present experimental results and analyze the performance of the proposed approach through a case study in space logistics. Finally, the paper is concluded in Section \ref{sec: Conclusion}, summarizing the findings and discussing potential future directions.

\section{Methodology}
\label{sec: Methodology}
\subsection{Mathematical Formulation}

In the realm of space logistics systems, MCNF describes space missions as multi-commodity flows inside a network, where nodes represent orbits, planets, or celestial objects, and arcs denote the trajectories linking the nodes. The spacecraft, crew, scientific instruments, and other payloads (e.g., water, food, oxygen, propellent) are all regarded as the commodities traversing the arcs \cite{chen_space_2019}.

The following shows the mathematical formulations for the time-expanded multi-commodity network flow problem, with the objective of optimizing the transportation of multiple commodities through a capacitated network. The network is represented as a directed graph \( G = (N, E) \), where \( N \) denotes the set of depots (nodes) and \( E \) represents the set of directed routes (edges) connecting the depots. Each commodity \( k \in K \) is transported across this network, with each commodity having a specific origin and destination depot. The supply and demand of these commodities vary over discrete time steps \( t \in T \). The transportation process involves vehicles with predefined capacity constraints, and the objective is to minimize the total transportation cost while satisfying flow conservation at each depot and time step and adhering to the vehicle capacity constraints on each route.

The flow of commodities is represented by \( x_{ij}^{k,t} \), which denotes the number of units of commodity \( k \) flowing along arc \( (i, j) \) at time \( t \). To support the flow, \( z_{ij}^t \) represents the number of vehicles deployed on route \( (i, j) \) at time \( t \). The cost of transportation is defined by \( c_{ij} \), which specifies the cost per vehicle for traversing arc \( (i, j) \). Each vehicle has a fixed capacity \( W \), and the load size of commodity \( k \) per unit flow is denoted by \( L_k \). The supply or demand of a commodity at a given depot \( i \) and time \( t \) is denoted as \( d_{ki}^t \). Additionally, a time delay \( \Delta t_{ij} \) is associated with traversing arc \( (i, j) \).
Based on the aforementioned notations, the time-expanded multi-commodity network flow problem is as follows: 

Minimize:
\begin{equation}
\mathcal{J}\ = \sum_{t} \sum_{(i,j) \in E} c_{ij} \cdot z_{ij}^t
\end{equation}

subject to the following constraints:

\begin{equation}
     \sum_{j \in N} x_{ij}^{k,t} \cdot L_k - \sum_{j \in N} x_{ji}^{k,t-\Delta t_{ji}} \cdot L_k = d_{ki}^t, \quad \forall i \in N, \, k \in K, \, t
\end{equation}

\begin{equation}
   \sum_{k \in K} x_{ij}^{k,t} \cdot L_k \leq z_{ij}^t \cdot W, \quad \forall (i,j) \in E, \, t
\end{equation}
\begin{equation}
x_{ij}^{k,t} \in \mathbb{Z}_{\geq 0}, \quad z_{ij}^t \in \mathbb{Z}_{\geq 0}, \quad \forall (i,j) \in E, \, k \in K, \, t
\end{equation}

The time-dependent supply or demand \( d_{ki}^t \) is defined as:
\[
d_{ki}^t =
\begin{cases} 
L_k, & \text{if } i \text{ is the origin depot for commodity } k \text{ at time } t, \\
-L_k, & \text{if } i \text{ is the destination depot for commodity } k \text{ at time } t, \\
0, & \text{otherwise (intermediate depots)}.
\end{cases}
\]

\subsection{Entropy Quantum Computing (EQC)}

Entropy quantum computing was introduced recently as an efficient optimization technique that leverages the principles of entropy minimization to guide quantum state evolution toward optimal solutions of a Hamiltonian \cite{nguyen_entropy_2024}. In a photonic realization, qudits are encoded as a superposition of photon numbers in time bins that evolve in an optical fiber loop which embeds a desired Hamiltonian as a dissipative operator. This effectively emulates imaginary time evolution, in which the propagation of higher-energy eigenstates is subject to dissipation and decoherence while the lower-energy eigenstates are promoted in this evolution. The current generation of the Entropy Quantum Computing (EQC) hardware at QCi, named Dirac-3, is a hybrid system that harnesses the high-speed parallel processing capabilities of photonics and the quantum nature of light with the precise control and programmability of electronic circuits. By leveraging photonic systems for entropy-driven state evolution—enabled by their natural ability to manipulate complex optical fields—and utilizing electronics for state initialization, feedback, and fine-tuning, hybrid entropy computing achieves efficient optimization. Such systems can explore large solution spaces rapidly through photonic modes propagation while employing electronic components to refine solutions iteratively. This synergy offers a pathway for solving challenging optimization problems in areas such as machine learning, signal processing, and network design, with the potential for high speed, low power consumption, and scalability.

The EQC hardware Dirac-3 targets optimization of a polynomial objective function $E$ in the following form:
\begin{equation}
    E = \sum_{i} C_i x_i +
        \sum_{i,j} J_{ij} x_i x_j +
        \sum_{i,j,k} T_{ijk} x_i x_j x_k +
        \sum_{i,j,k,l} Q_{ijkl} x_i x_j x_k x_l +
        \sum_{i,j,k,l,m} P_{ijklm} x_i x_j x_k x_l x_m .
\end{equation}
where, $x_i$ are optimization variables, $C_i$ are real-valued coefficients of linear terms, and $J_{ij}$, $T_{ijk}$, $Q_{ijkl}$, $P_{ijklm}$ respectively represent two-body to fifth-body interaction coefficients that are real numbers subject to the tensors $J$, $T$, $Q$, and $P$ being symmetric under all permutations of the indices. Governed by the physical mechanism of operation of the EQC hardware, $x_i$ represents photon numbers in the $i$’th time bin degrees of freedom. Thus, the optimization variables are non-negative real numbers, $x_i \geq 0$, and are subject to a sum constraint $\sum_i x_i = R$.

A major advantage of EQC compared to other quantum annealers is its flexible variable encoding capability \cite{nguyen_entropy_2024}. While, in the form discussed above, EQC allows for solving continuous-variable optimization problems, by allocating multiple photon time bins to one variable, one can effectively encode binary or integer variables for solving combinatorial optimization problems in the forms of binary, integer or even mixed-integer problems.

\section{Results and Analysis}
\label{sec: Results and Analysis}

\subsection{Problem Settings}

\begin{figure}[H]
    \centering
    \includegraphics[width=0.5\linewidth]{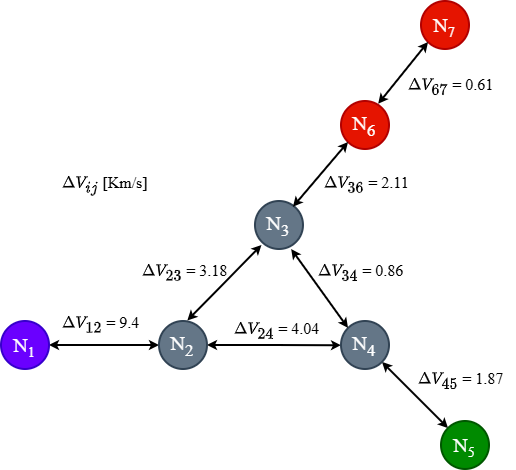}
    \caption{Network model}
    \label{fig:network}
\end{figure}

The network presented in Fig \ref{fig:network} highlights a hypothetical scenario involving orbital locations such as Earth (\( N_1 \)), Low Earth Orbit (LEO) (\( N_2 \)), Lunar Transfer Orbit (LTO) (\( N_3 \)), Low Lunar Orbit (LLO) (\( N_4 \)), Lunar Surface (LS) (\( N_5 \)), Low Mars Orbit (LMO) (\( N_6 \)), and Mars (\( N_7 \)). Each arc's \( \Delta V_{ij} \) is labeled on the figure, with values indicating the velocity changes required for traversing specific trajectories, measured in kilometers per second (\( \text{Km/s} \)). In this problem, the cost-efficient parameter \( c_{ij} \) for each arc \( (i, j) \) is equivalent to the (\( \Delta V_{ij} \)), representing the velocity change required to traverse the arc. The time of flight (\( \Delta t_{ij} \)) is fixed and equal to \( 1 \) for all arcs, simplifying the formulation and ensuring consistency across the network.

\begin{equation}
c_{ij} \sim \Delta V_{ij}, \quad \forall (i, j) \in E
\end{equation}

The roles of the nodes, the commodities being transported, and the time-dependent supply and demand (\( d_{ki}^t \)) for each commodity across discrete time steps are outlined in Table~\ref{tab:nodes_table}. Each node is assigned a specific role—origin, intermediate, or destination—and supports the transportation of two commodity types, \( L_1 \) and \( L_2 \). The table provides a structured view of how supply and demand evolve across the network over time. Also, the essential parameters used in the problem formulation are presented in Table~\ref{tab:parameters_table}. These include the cost parameter \( c_{ij} \), vehicle capacity \( W \), and load sizes \( L_1 \) and \( L_2 \) for each commodity type. Together, these parameters establish the framework for modeling and solving the multicommodity flow problem in a space logistics context.

\begin{table}[h!]
\centering
\caption{Nodes Table}
\label{tab:nodes_table}
\begin{tabular}{llcl}
\toprule
\textbf{Node} & \textbf{Role} & \textbf{Commodity} & \textbf{Supply/Demand (\(d_{i,k}^t\)) for \(t = 1,2,3,4,5,6\)} \\ 
\midrule
\(N_1\) (Earth) & Origin & \((L_1, L_2)\) & \( (+40,+80), (+60,+120), (0,0), (0,0), (0,0), (0,0) \) \\
\(N_2\) (LEO)   & Intermediate & \((L_1, L_2)\) & \( (0,0), (0,0), (0,0), (0,0), (0,0), (0,0) \) \\
\(N_3\) (LTO)   & Intermediate & \((L_1, L_2)\) & \( (0,0), (0,0), (0,0), (0,0), (0,0), (0,0) \) \\
\(N_4\) (LLO)   & Intermediate & \((L_1, L_2)\) & \( (0,0), (0,0), (0,0), (0,0), (0,0), (0,0) \) \\
\(N_5\) (LS)    & Destination & \((L_1, L_2)\) & \( (0,0), (0,0), (0,0), (0,0), (-20,-40), (-30,-60) \) \\
\(N_6\) (LMO)   & Intermediate & \((L_1, L_2)\) & \( (0,0), (0,0), (0,0), (0,0), (0,0), (0,0) \) \\
\(N_7\) (Mars)  & Destination & \((L_1, L_2)\) & \( (0,0), (0,0), (0,0), (0,0), (-20,-40), (-30,-60) \) \\
\bottomrule
\end{tabular}
\end{table}

\newcolumntype{L}[1]{>{\raggedright\arraybackslash}p{#1}}

\begin{table}[h!]
\centering
\caption{Parameters Table}
\label{tab:parameters_table}
\begin{tabular}{@{}L{3cm}L{5.5cm}L{4cm}@{}}
\toprule
\textbf{Parameter} & \textbf{Description} & \textbf{Value} \\ 
\midrule
\( c_{ij} \)       & Cost for each arc, defined as \( \Delta V_{ij} \) & Listed in the arcs table \\ 
\( \Delta t_{ij} \) & Time of flight for each arc                      & \( 1 \)                  \\ 
\( W \)            & Vehicle capacity (units)                         & \( 100 \)                \\ 
\( L_1, L_2 \)     & Load size per unit flow (units)                  & \( L_1 = 10, L_2 = 20 \) \\ 
\bottomrule
\end{tabular}
\end{table}
\FloatBarrier

\subsection{Quantum Results}

\input{quantumresults}

\section{Conclusion}
\label{sec: Conclusion}
In this paper, we present the potential of quantum computing, specifically entropy quantum computing, in addressing the complex challenges of space logistics mission planning. By formulating the problem as a time-dependent multicommodity network flow and leveraging quantum Hamiltonian models, we efficiently explored the vast solution space required for interplanetary missions. The use of entropy quantum computing allowed for the simultaneous evaluation of multiple solutions, with the potential to significantly reduce computation time while adhering to vehicle capacity, resource, and temporal constraints.
The results validate the feasibility of quantum-assisted frameworks in optimizing logistical networks that are intractable for classical methods. The approach successfully planned routes, managed resource allocations, and ensured the timely delivery of commodities across a dynamic space logistics network. These findings underline the transformative capabilities of quantum computing in revolutionizing mission planning for complex aerospace systems. Future work will focus on enhancing the scalability of the approach with advanced quantum algorithms and hardware, as well as exploring its applicability to broader domains, such as terrestrial logistics and supply chain optimization.



\section*{Acknowledgments}

This material is based upon work supported by the National Science Foundation under Grant No. 2301627. Any opinions, findings, and conclusions or recommendations expressed in this material are those of the authors and do not necessarily reflect the views of the National Science Foundation.

 \newpage
\bibliography{references}

\end{document}

%% file: quantumresults.tex
We have developed a method using the Python programming language to accept instances of logistics models as described in this paper and create a Hamiltonian operator whose ground state describes the optimal solution to the problem. This code formulates the objective function as a vector, $P_1$, and all the constraints from the model description as a linear system of equalities in left-hand and right-hand side pairs of arrays. These constraints become penalties, $P_2$, a polynomial of order 2 just as is done for QUBO in \cite{glover_tutorial_2019}. The Hamiltonian is described as
\begin{equation}
    H(x,z)=P_1(x,z)+\alpha P_2(x,z)
\end{equation}
where $\alpha$ is a multiplier that is chosen to elevate the penalties enough to make positive penalties greater than the objective function coefficients. Increasing this value provides more separation of feasible solutions from infeasible ones in terms of energy value, but this also can reduce the impact of small coefficients in the objective function. The ratio of the magnitudes of the largest to the smallest absolute values in the problem is called the dynamic range of the problem and is reported in decibels. This means that if the ratio is 200, then the dynamic range is $\log_{10} 200\approx 23$dB. In total, $H(x,z)$ has 116 variables, 89 of which are decision variables and 29 are slack variables. Among all variables, there are 331 levels, each a time dimension in Dirac-3. The coefficients have a dynamic range of 30.78 dB.

In addition to selecting the $\alpha$ parameter for the model formulation, Dirac-3 provides configuration parameters for manipulating the sampler behavior. The parameters used in this example are in table \ref{tab:dirac-3-params}.

\begin{table}[h!]
\centering
\caption{Device Configuration Parameters}
\label{tab:dirac-3-params}
\begin{tabular}{@{}L{7cm}lll@{}L{3cm}}
\toprule
\textbf{Name}                            & \textbf{Value} \\ 
\midrule
Quantum Fluctuation Coefficient                               & $1/\sqrt{7}$              \\ 
Relaxation Schedule                                            & 2              \\ 
\bottomrule
\end{tabular}
\end{table}

\begin{figure}[h!]
    \centering
    \begin{minipage}[t]{0.49\textwidth}
        \centering
        \includegraphics[width=\linewidth]{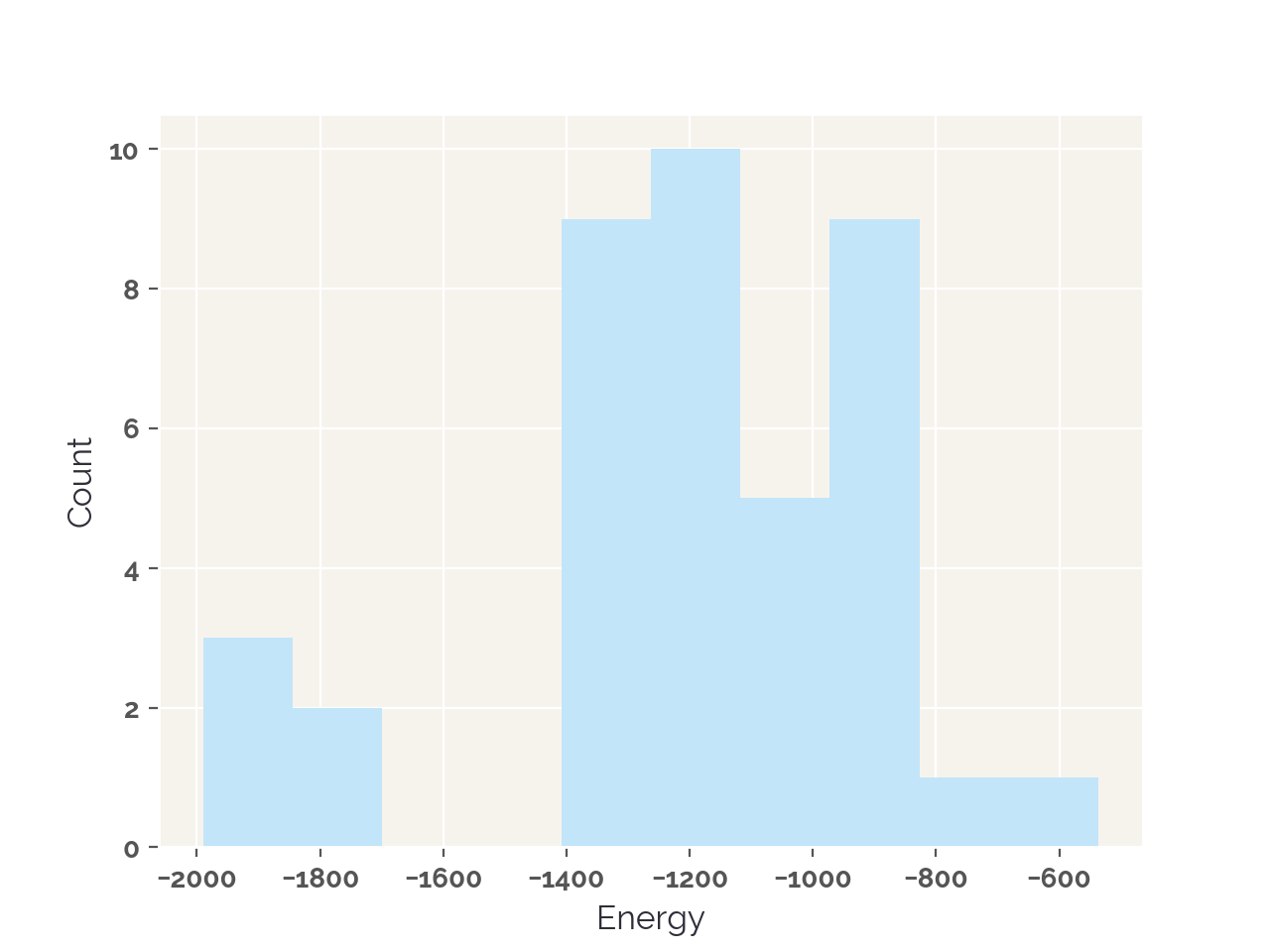}
        \caption{Distribution of Sample Energies}
        \label{fig:distribution}
    \end{minipage}
    \hfill
    \begin{minipage}[t]{0.49\textwidth}
        \centering
        \includegraphics[width=\linewidth]{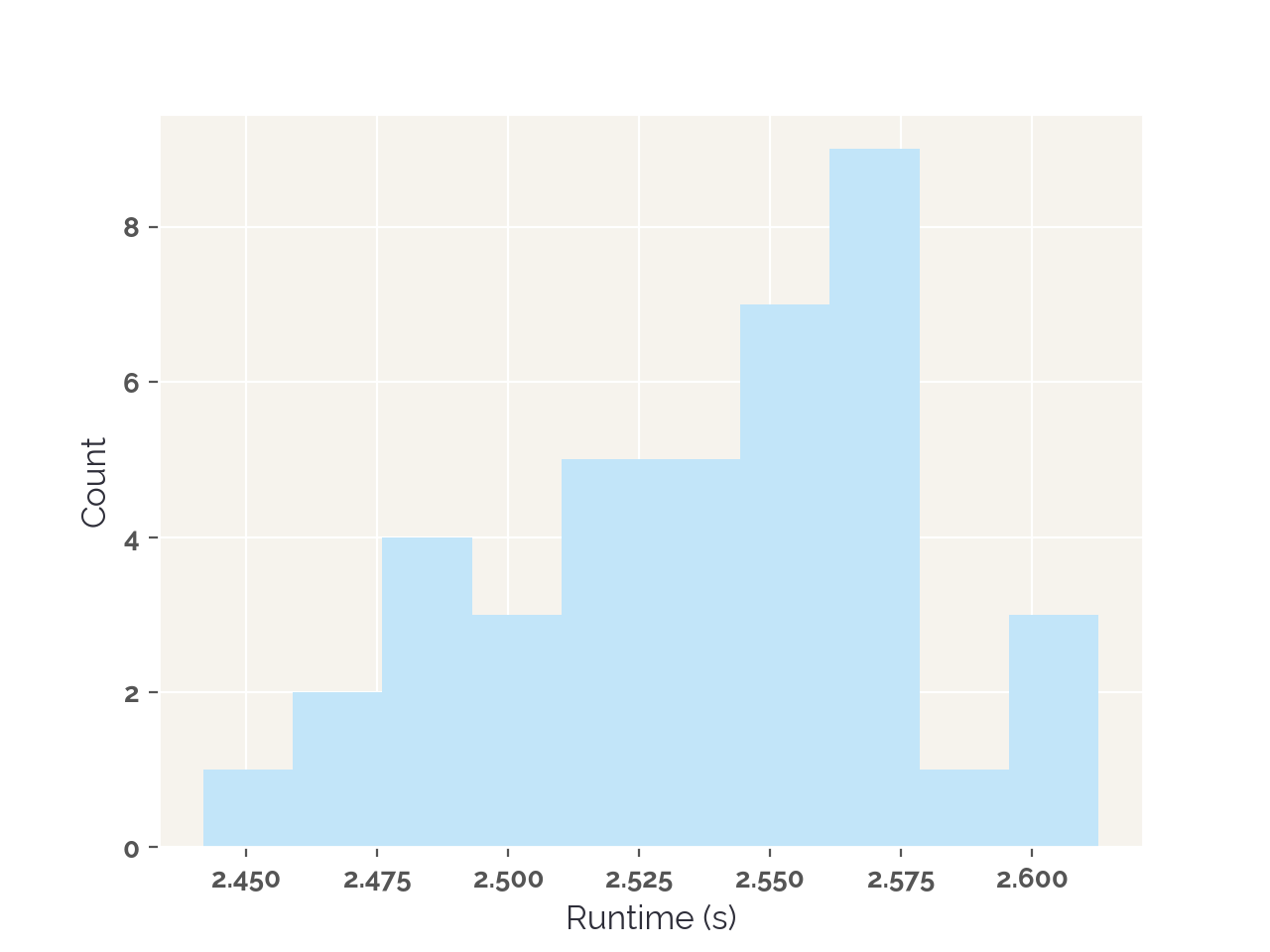}
        \caption{Distribution of Sample Runtimes}
        \label{fig:times}
    \end{minipage}
\end{figure}

The computational results were obtained from Dirac-3 via the Cloud application programming interface. The solution displayed in Tables \ref{tab:flow-table}, \ref{tab:vehicle-table}, and \ref{tab:total-cargo} is the best energy (-1989.2251) solution out of 40 samples taken using the parameters listed. See Figure \ref{fig:distribution} for the distribution of energies for these 40 solutions. Each of the samples took between 2.44 and 2.61 seconds of device and post processing time.  The average time was 2.53 seconds as seen in Figure \ref{fig:times}. The device time was approximately 1.5 seconds per sample and post processing was approximately 1 second per sample. Post processing time is time for auxilary computations that occur during each call to Dirac hardware including intermediate energy calculations and objective function adjustments. The commodity flow table shows the quantity of a commodity that is at each node at the end of each time step. For clarity, the demand quantities are not removed from the nodes (\(N_5\) and \(N_7\)), so at $t=6$, the total delivery volume is displayed. The vehicle movement table lists the count of vehicles traversing an arc at each time step. Together with the total cargo movements table, the vehicle solution is seen to be consistent with the cargo size. The raw solution contained a single load of \(L_1\) from \(N_1\) to \(N_2\), \(N_2\) to \(N_3\), and \(N_3\) to \(N_4\) in time steps 4, 5, and 6. It also contained a single load movement of \(L_2\) from \(N_2\) to \(N_3\) in time step 5. This prompted the inclusion of a post-processing method which set all commodity flows to 0 when no vehicle was scheduled (negligible to runtime). This eliminated the spurious movements of commodities with no supply or demand, while not modifying the vehicle schedule determined by the device. The resulting solution satisfies all demand with supply movements at appropriate time steps with an objective function value of 62.08 Km/s. This is not the optimal solution because there is clearly one more vechicle than is required at time step 4 from \(N_3\) to \(N_4\). As the cost from \(N_3\) to \(N_4\) is one of the smallest coefficients of the model, including the two cost coefficients less than or equal to 0.86 (\(N_3\)-\(N_4\) and \(N_6\)-\(N_7\)) increase the dynamic range of the Hamiltonian to 30.78 db, which approaches close to the range limit of the Dirac-3 machine. The Hamiltonian without these small coeffients has a dynamic range of 25.9 dB.

\begin{table}[H]
\centering
\begin{minipage}[t]{0.48\textwidth} 
\centering
\caption{Vehicle-Arc Traversals by Time Step}

\begin{tabular}{@{}lcccccc@{}}
\toprule
\multirow{2}{*}{\textbf{Arc}} & \multicolumn{6}{c}{\textbf{Time}} \\ 
\cmidrule(lr){2-7}
  & \textbf{1} & \textbf{2} & \textbf{3} & \textbf{4} & \textbf{5} & \textbf{6} \\ 
\midrule
\(N_1\), \(N_2\) & 0 & 2 & 2 & 0 & 0 & 0 \\
\(N_2\), \(N_3\) & 0 & 0 & 2 & 2 & 0 & 0 \\
\(N_2\), \(N_4\) & 0 & 0 & 0 & 0 & 0 & 0 \\
\(N_3\), \(N_6\) & 0 & 0 & 0 & 1 & 1 & 0 \\
\(N_3\), \(N_4\) & 0 & 0 & 0 & 2 & 1 & 0 \\
\(N_4\), \(N_5\) & 0 & 0 & 0 & 0 & 1 & 1 \\
\(N_6\), \(N_7\) & 0 & 0 & 0 & 0 & 1 & 1 \\
\(N_4\), \(N_3\) & 0 & 0 & 0 & 0 & 0 & 0 \\
\bottomrule
\label{tab:vehicle-table}
\end{tabular}
\end{minipage}
\hfill
\begin{minipage}[t]{0.48\textwidth} 
\centering
\caption{Total Cargo Movements}
\label{tab:total-cargo}
\begin{tabular}{@{}lcccccc@{}}
\toprule
\multirow{2}{*}{\textbf{Arc}} & \multicolumn{6}{c}{\textbf{Time}} \\ 
\cmidrule(lr){2-7}
  & \textbf{1} & \textbf{2} & \textbf{3} & \textbf{4} & \textbf{5} & \textbf{6} \\ 
\midrule
\(N_1\), \(N_2\) & 0 & 120 & 180 & 0 & 0 & 0 \\
\(N_2\), \(N_3\) & 0 & 0 & 120 & 180 & 0 & 0 \\
\(N_2\), \(N_4\) & 0 & 0 & 0 & 0 & 0 & 0 \\
\(N_3\), \(N_6\) & 0 & 0 & 0 & 0 & 60 & 90 \\
\(N_3\), \(N_4\) & 0 & 0 & 0 & 0 & 60 & 90 \\
\(N_4\), \(N_5\) & 0 & 0 & 0 & 0 & 60 & 90 \\
\(N_6\), \(N_7\) & 0 & 0 & 0 & 0 & 60 & 90 \\
\(N_4\), \(N_3\) & 0 & 0 & 0 & 0 & 0 & 0 \\
\bottomrule
\end{tabular}
\end{minipage}
\end{table}

\begin{table}[]
\centering
\caption{Commodity Flow By Time Step (Demand Remains In Inventory)}
\label{tab:flow-table}
\begin{tabular}{@{}ccccccc@{}}
\toprule
\multicolumn{1}{c}{\multirow{2}{*}{\textbf{Node, Commodity}}} & \multicolumn{6}{c}{\textbf{Time}}                                           \\ \cmidrule(l){2-7} 
\multicolumn{1}{c}{}                                          & \textbf{1} & \textbf{2} & \textbf{3} & \textbf{4} & \textbf{5} & \textbf{6} \\ \midrule
\(N_1\), \(L_1\) & 40 & 60 & 0  & 0  & 0  & 0  \\
\(N_2\), \(L_1\) & 0  & 40 & 60 & 0  & 0  & 0  \\
\(N_3\), \(L_1\) & 0  & 0  & 40 & 60 & 0  & 0  \\
\(N_4\), \(L_1\) & 0  & 0  & 0  & 20 & 30 & 0  \\
\(N_5\), \(L_1\) & 0  & 0  & 0  & 0  & 20 & 50 \\
\(N_6\), \(L_1\) & 0  & 0  & 0  & 20 & 30 & 0  \\
\(N_7\), \(L_1\) & 0  & 0  & 0  & 0  & 20 & 50 \\
\(N_1\), \(L_2\) & 80 & 120 & 0  & 0  & 0  & 0  \\
\(N_2\), \(L_2\) & 0  & 80  & 120 & 0  & 0  & 0  \\
\(N_3\), \(L_2\) & 0  & 0   & 80  & 120 & 0  & 0  \\
\(N_4\), \(L_2\) & 0  & 0   & 0   & 40  & 60  & 0  \\
\(N_5\), \(L_2\) & 0  & 0   & 0   & 0   & 40  & 100 \\
\(N_6\), \(L_2\) & 0  & 0   & 0   & 40  & 60  & 0  \\
\(N_7\), \(L_2\) & 0  & 0   & 0   & 0   & 40  & 100 \\
\bottomrule
\label{tab:flow-table}
\end{tabular}
\end{table}